\documentclass[12pt]{amsart}

\usepackage{amssymb, amscd, txfonts}
\usepackage{graphicx}

\numberwithin{equation}{section}

\newtheorem{theorem}{Theorem}[section]
\newtheorem{proposition}[theorem]{Proposition}
\newtheorem{lemma}[theorem]{Lemma}

\theoremstyle{definition}

\theoremstyle{remark}

\newcommand{\C}{\mathbb{C}}
\newcommand{\F}{\mathbb{F}}
\newcommand{\proj}{{\mathbb P}}

\newcommand{\PGL}{{\rm PGL}}
\newcommand{\Oline}{\mathcal{O}_{{\mathbb P}^{1}}}
\newcommand{\Oplane}{\mathcal{O}_{{\mathbb P}^{2}}}
\newcommand{\Orel}{\mathcal{O}_{\pi}}

\newcommand{\autF}{{\aut}(\mathbb{F}_1)}

\DeclareMathOperator{\aut}{Aut}

\begin{document}

\title[]{The rationality of the moduli spaces of trigonal curves of odd genus}
\author[]{Shouhei Ma}
\address{Graduate~School~of~Mathematical~Sciences, the~University~of~Tokyo, 3-8-1~Komaba, Meguro-ku, Tokyo 153-8914, Japan}
\email{sma@ms.u-tokyo.ac.jp}
\maketitle 

\begin{abstract}
The moduli spaces of trigonal curves of odd genus $g\geq5$ are proven to be rational. 
\end{abstract}

\maketitle


\section{Introduction}\label{sec: intro}

The object of this article is to prove the following. 

\begin{theorem}\label{main}
The moduli space $\mathcal{T}_g$ of trigonal curves of genus $g=2n+1$ with $n\geq2$ is rational. 
\end{theorem}

By a \textit{trigonal curve} we mean an irreducible smooth projective curve which admits a degree $3$ morphism to ${\proj}^1$. 
A trigonal curve of genus $g\geq5$ has a unique $g_3^1$, so that 
the space $\mathcal{T}_g$ to be studied is regarded as a sublocus of $\mathcal{M}_g$, 
the moduli space of curves of genus $g$. 
Shepherd-Barron \cite{SB} proved the rationality of $\mathcal{T}_g$ for $g=4n+2$ with $n\geq1$. 
Hence the space $\mathcal{T}_g$ is rational possibly except when the genus $g$ is divisible by $4$. 
For the one lower gonality, 
Katsylo and Bogomolov \cite{Ka1}, \cite{B-K} established the rationality of the moduli spaces of hyperelliptic curves.

The proof of Theorem \ref{main} is based on the classical relation between trigonal curves and 
the Hirzebruch surfaces ${\F}_N={\proj}({\Oline}\oplus{\Oline}(N))$. 
Recall that a canonically embedded trigonal curve $C\subset{\proj}^{g-1}$ of genus $g\geq5$ 
lies on a unique rational normal scroll $S$. 
The scroll $S$ may obtained either as the intersection of quadrics containing $C$, 
or as the scroll swept out by the lines spanned by the fibers of the trigonal map. 
The surface $S$ is the image of a Hirzebruch surface ${\F}_N$ 
by a linear system $|{\Orel}(1)\otimes\pi^{\ast}{\Oline}(a)|$, $a>0$, 
where $\pi:{\F}_N\to{\proj}^1$ is the natural projection. 
The trigonal map of $C$ is the restriction of $\pi$. 
When $C$ is general in the moduli $\mathcal{T}_g$, 
we have $N=0$ or $1$ depending on whether $g$ is even or odd. 
Thus, if $L_{3,b}$ denotes the line bundle ${\Orel}(3)\otimes\pi^{\ast}{\Oline}(b)$ on ${\F}_1$ 
with $g=2b+1$, we have the birational equivalence 
\begin{equation}
\mathcal{T}_g \sim |L_{3,b}|/{\aut}({\F}_1). 
\end{equation}
Here $|L_{3,b}|/{\aut}({\F}_1)$ stands for a rational quotient of the linear system $|L_{3,b}|$ by the algebraic group ${\aut}({\F}_1)$. 
Then Theorem \ref{main} is equivalent to the following assertion in invariant theory. 

\begin{theorem}\label{main'}
For the line bundle $L_{3,b}$ on the Hirzebruch surface ${\F}_1$ 
the quotient $|L_{3,b}|/{\aut}({\F}_1)$ is rational for $b\geq2$. 
\end{theorem}

The rest of this article is devoted to the proof of this theorem. 
In Section \ref{sec: symmetric product} we construct an ${\autF}$-equivariant map from $|L_{3,b}|$ to $S^b{\F}_1$, 
the symmetric product of ${\F}_1$, which plays crucial role in the proof. 
In Section \ref{sec: g>7} the rationality for $g\geq9$ is established by using the rational normal curves. 
In Section \ref{sec: g<9} the rationality of $\mathcal{T}_7$ and $\mathcal{T}_5$ is proved,

Throughout this article we work over the field of complex numbers. 
We denote by $\pi:{\F}_1\to{\proj}^1$ the natural projection. 
The $(-1)$-curve on ${\F}_1$ is denoted by $\Sigma$.  
The line bundle ${\Orel}(a)\otimes\pi^{\ast}{\Oline}(b)$ on ${\F}_1$ will be written as $L_{a,b}$. 
The bundle ${\Orel}(1)$ is the pullback of ${\Oplane}(1)$ by the blow-down ${\F}_1\to{\proj}^2$.


\section{Symmetric product of the Hirzebruch surface}\label{sec: symmetric product}

Let ${\proj}\mathcal{E}$ be the projective space bundle ${\proj}\pi_{\ast}{\Orel}(2)$ on ${\proj}^1$. 
The variety ${\proj}\mathcal{E}$ parametrizes unordered pairs $q_++q_-$ of two points of ${\F}_1$ 
which lie on the same $\pi$-fiber. 
We have a rational map 
\begin{equation}\label{map 1}
\varphi_1 : |L_{3,b}| \dashrightarrow S^b({\proj}\mathcal{E}), \quad  C\mapsto \sum_{i=1}^b(q_{i+}+q_{i-}) 
\end{equation}
defined as follows. 
If $C|_{\Sigma}=p_1+\cdots+p_b$ and $F_i$ is the $\pi$-fiber passing $p_i$, 
we set $q_{i+}+q_{i-}=C|_{F_i}-p_i$. 
The map $\varphi_1$ is clearly ${\autF}$-equivariant. 
Next we define a rational map 
\begin{equation}\label{map 2}
\varphi_2 : S^b({\proj}\mathcal{E}) \dashrightarrow S^b{\F}_1, \quad  \sum_{i=1}^b(q_{i+}+q_{i-})\mapsto \sum_{i=1}^{b}q_i 
\end{equation}
as follows. 
If $F_i$ is the $\pi$-fiber passing $\{ q_{i+}, q_{i-} \}$ and $p_i=F_i\cap\Sigma$, 
there exists a unique involution $\iota_i$ of $F_i\simeq{\proj}^1$ which fixes $p_i$ and interchanges $q_{i+}$ and $q_{i-}$. 
Then we let $q_i\in F_i$ be the fixed point of $\iota_i$ other than $p_i$. 
By the uniqueness of $\iota_i$ the map $\varphi_2$ is ${\autF}$-equivariant. 
We study the composition map 
\begin{equation}\label{map 3}
\varphi = \varphi_2 \circ \varphi_1 : |L_{3,b}| \dashrightarrow S^b{\F}_1. 
\end{equation}

\begin{lemma}\label{lemma 2.1}
The map $\varphi$ is dominant with a general fiber being an open set of a linear subspace of $|L_{3,b}|$. 
\end{lemma}

\begin{proof}
For a general point $q_1+\cdots+q_b \in S^b{\F}_1$ let $F_i$ be the $\pi$-fiber passing $q_i$ and let $p_i=F_i\cap\Sigma$. 
We take an inhomogeneous coordinate $x_i$ of $F_i\simeq{\proj}^1$ 
in which $p_i$ is $\{ x_i=0 \}$ and $q_i$ is $\{ x_i=\infty \}$. 
The involution of $F_i$ fixing $p_i$ and $q_i$ is given by $x_i\mapsto -x_i$. 
A smooth curve $C\in |L_{3,b}|$ is contained in $\varphi^{-1}(q_1+\cdots+q_b)$ if and only if 
$C|_{F_i}$ has the equation $x_i(\alpha_ix_i^2+\beta_i)=0$ for each $i=1,\cdots, b$. 
Since these are $2b$ linear conditions on $|L_{3,b}|$, 
namely the vanishing of the coefficient of $x_i^2$ and the constant term for $C|_{F_i}$, 
the second assertion is proved. 
The dominancy of $\varphi$ is a consequence of the dimension counting 
${\dim}|L_{3,b}|=4b+9>2b$. 
\end{proof}

\begin{lemma}\label{lemma 2.2}
The group ${\autF}$ acts on $S^b{\F}_1$ almost freely if $b\geq4$. 
\end{lemma}

\begin{proof}
First we treat the case $b\geq5$. 
If a general point $p_1+\cdots+p_b \in S^b{\F}_1$ is fixed by a $g\in{\autF}$, 
then $g$ stabilizes a general $b\geq5$ point set of the $(-1)$-curve $\Sigma$ so that $g$ acts trivially on $\Sigma$. 
Hence $g$ fixes each $p_i$. 
As ${\autF}$ acts almost freely on $({\F}_1)^b$, it follows that $g={\rm id}$. 

Next we study the case $b=4$. 
Let $f:{\F}_1\to{\proj}^2$ be the blow-down. 
For a general $p_1+\cdots+p_4\in S^4{\F}_1$ 
there exists a unique smooth conic $Q$ passing $f(\Sigma)$ and $f(p_1), \cdots, f(p_4)$. 
Any $g\in{\autF}$ fixing $p_1+\cdots+p_4$, 
regarded as an element of ${\PGL}_3$, 
preserves $Q$ and the five point set $f(\Sigma), f(p_1), \cdots, f(p_4)$ on it. 
Hence $g$ acts trivially on $Q$, which implies that $g={\rm id}$. 
\end{proof}

We shall apply the no-name lemma (see \cite{Do}, and also \cite{C-G-R} for non-reductive groups) 
to the map $\varphi$ when $b\geq4$. 
For that we note the following. 

\begin{lemma}\label{linearization} 
Every line bundle on ${\F}_1$ admits an ${\autF}$-linearization. 
\end{lemma}

\begin{proof}
We have canonical ${\autF}$-linearizations on the bundles 
$K_{{\F}_1}=L_{-2,-1}$, $\pi^{\ast}K_{{\proj}^1}=L_{0,-2}$, and $f^{\ast}K_{{\proj}^2}=L_{-3,0}$ 
where $f:{\F}_1\to{\proj}^2$ is the blow-down of $\Sigma$. 
These induce ${\autF}$-linearizations on $L_{1,0}$ and $L_{0,1}$. 
Since ${\rm Pic}({\F}_1)$ is freely generated by $L_{1,0}$ and $L_{0,1}$, the lemma is proved. 
\end{proof}

By Lemma \ref{linearization} 
the ${\autF}$-action on $|L_{3,b}|$ is induced by an ${\autF}$-representation on $H^0(L_{3,b})$. 
Then Lemma \ref{lemma 2.1} shows that $|L_{3,b}|$ is ${\autF}$-birational to 
the projectivization of an ${\autF}$-linearized vector bundle on an open set of $S^b{\F}_1$. 
By Lemma \ref{lemma 2.2} we may apply the no-name lemma to see the 

\begin{proposition}\label{prop 2.4}
For $b\geq4$ we have a birational equivalence 
\begin{equation}
|L_{3,b}|/{\autF} \sim {\proj}^{2b+9}\times(S^b{\F}_1/{\autF}). 
\end{equation}
\end{proposition}

Thus the rationality of $|L_{3,b}|/{\autF}$ for $b\geq4$ is reduced to a stable rationality of $S^b{\F}_1/{\autF}$.


\section{Projection of rational normal curve}\label{sec: g>7}

In this section we prove a stable rationality of the quotient $S^b{\F}_1/{\autF}$ 
to derive Theorem \ref{main'} for $b\geq4$. 
For an integer $d\geq0$ we consider the universal curve $f:\mathcal{H}_d \to |L_{1,d}|$ over the linear system $|L_{1,d}|$. 
The variety $\mathcal{H}_d$ is defined as a divisor on ${\F}_1\times|L_{1,d}|$, 
and $f$ is the restriction of the second projection ${\F}_1\times|L_{1,d}| \to |L_{1,d}|$. 
The bundle $L_{0,1}$ on ${\F}_1$ induces a relative hyperplane bundle for $f$ which we denote by 
$\mathcal{O}_f(1)$. 
Let 
\begin{equation}
\mathcal{H}_{d,b}={\proj}f_{\ast}\mathcal{O}_f(b). 
\end{equation}
An open set of $\mathcal{H}_{d,b}$ parametrizes pairs $(H, q_1+\cdots+q_b)$ 
where $H\in|L_{1,d}|$ is smooth and $q_1, \cdots, q_b$ are $b$ points on $H$. 
Note that a smooth $H\in|L_{1,d}|$ is a section of $\pi$.

\begin{lemma}\label{lemma 3.1}
For $4\leq b\leq2d+2$ we have a birational equivalence 
\begin{equation}\label{eqn 3.1}
\mathcal{H}_{d,b}/{\autF} \sim {\proj}^{2d+2-b}\times(S^b{\F}_1/{\autF}). 
\end{equation}
\end{lemma}

\begin{proof}
Consider the evaluation map 
\begin{equation}
\psi : \mathcal{H}_{d,b} \dashrightarrow S^b{\F}_1, \quad  (H, q_1+\cdots+q_b) \mapsto q_1+\cdots+q_b. 
\end{equation}
The fiber $\psi^{-1}(q_1+\cdots+q_b)$ over a general $q_1+\cdots+q_b$ is an open set of the sub linear system of 
$|L_{1,d}|$ of curves passing $q_1, \cdots, q_b$. 
Since ${\dim}|L_{1,d}|=2d+2\geq b$, 
$\psi^{-1}(q_1+\cdots+q_b)$ is non-empty and of dimension $2d+2-b$. 
In particular, $\psi$ is dominant. 
Then we may apply the no-name lemma for $\psi$ as like the proof of Proposition \ref{prop 2.4} 
to deduce the equivalence \eqref{eqn 3.1}. 
\end{proof}

By a comparison of Proposition \ref{prop 2.4} and Lemma \ref{lemma 3.1}, 
it suffices for the proof of Theorem \ref{main'} for $b\geq4$ to show the rationality of $\mathcal{H}_{d,b}/{\autF}$ 
for \textit{one} $d$ in the range $b\leq 2d+2\leq 3b+9$. 
We begin with the 

\begin{lemma}\label{lemma 3.2}
For $d\geq5$ we have a birational equivalence 
\begin{equation}
\mathcal{H}_{d,b}/{\autF} \sim {\proj}^b\times(|L_{1,d}|/{\autF}). 
\end{equation}
\end{lemma}

\begin{proof}
This lemma is an application of the no-name method for the fibration $\mathcal{H}_{d,b} \to |L_{1,d}|$. 
Since the bundle $L_{0,1}$ on ${\F}_1$ admits an ${\autF}$-linearization, 
so is the bundle $\mathcal{O}_f(1)$ on the universal curve $\mathcal{H}_d$. 
Hence the sheaf $f_{\ast}\mathcal{O}_f(b)$ on $|L_{1,d}|$ is ${\autF}$-linearized. 
It remains to check the almost freeness of the ${\autF}$-action on $|L_{1,d}|$ for $d\geq5$. 
For a general $H\in|L_{1,d}|$ the intersection $H\cap\Sigma$ is a general $d$ point set of $H\simeq {\proj}^1$. 
If a $g\in {\autF}$ stabilizes $H$, then we have $g(H\cap\Sigma)=H\cap\Sigma$ so that $g$ acts trivially on $H$. 
This is enough for concluding that $g={\rm id}$. 
\end{proof}

Blowing-down ${\F}_1$ to ${\proj}^2$, 
we see that the quotient $|L_{1,d}|/{\autF}$ is birational to the ${\PGL}_3$-quotient of 
the space $\mathcal{X}_d$ of rational plane curves of degree $d+1$ having an ordinary $d$-fold point. 
Let $\widetilde{\mathcal{X}}_d$ be the space of morphisms $\phi:{\proj}^1\to{\proj}^2$ such that 
$\phi^{\ast}{\Oplane}(1)\simeq{\Oline}(d+1)$ and $\phi({\proj}^1)\in \mathcal{X}_d$. 
We have 
\begin{equation}
|L_{1,d}|/{\autF} \sim {\PGL}_2\backslash\widetilde{\mathcal{X}}_d/{\PGL}_3. 
\end{equation}

Let ${\proj}V_{d+1}=|{\Oline}(d+1)|^{\vee}$ and 
$\Gamma_{d+1}\subset {\proj}V_{d+1}$ be the rational normal curve $\phi_0({\proj}^1)$ 
where $\phi_0$ is the embedding associated to ${\Oline}(d+1)$. 
Recall that every morphism $\phi:{\proj}^1\to{\proj}^2$ with $\phi^{\ast}{\Oplane}(1)\simeq{\Oline}(d+1)$ 
is the composition of 
$(1)$ the isomorphism $\phi_0:{\proj}^1 \to \Gamma_{d+1}$, 
$(2)$ the projection $\Gamma_{d+1}\to {\proj}(V_{d+1}/W)$ from a $(d-2)$-plane ${\proj}W\subset {\proj}V_{d+1}$ 
         which is disjoint from $\Gamma_{d+1}$, and 
$(3)$ an isomorphism ${\proj}(V_{d+1}/W)\to {\proj}^2$. 
The group ${\PGL}_3$ acts on $\widetilde{\mathcal{X}}_d$ by transformation of an isomorphism ${\proj}(V_{d+1}/W)\to {\proj}^2$. 
Hence the quotient $\widetilde{\mathcal{X}}_d/{\PGL}_3$ is naturally birational to the locus 
$\mathcal{Y}_d$ in the Grassmannian  $\mathbf{G}(d-2, {\proj}V_{d+1})$ consisting 
of $(d-2)$-planes ${\proj}W$ such that 
(i) ${\proj}W\cap\Gamma_{d+1}=\emptyset$ and 
(ii) there exists a $(d-1)$-plane ${\proj}U$ containing ${\proj}W$ with ${\proj}U\cap\Gamma_{d+1}$ being a $d$ point set. 
For such a ${\proj}W$, the $(d-1)$-plane ${\proj}U$ is spanned by the point set ${\proj}U\cap\Gamma_{d+1}$ 
because of the fact that any distinct $d$ points on a rational normal curve in ${\proj}^{d+1}$ are  linearly independent. 
Also ${\proj}U$ is uniquely determined by ${\proj}W$ for an irreducible plane curve of degree $d+1$ has 
at most one singularity of multiplicity $d$. 
These two facts imply that $\mathcal{Y}_d$ is identified with an open set of the locus 
\begin{equation}
\mathcal{Z}_d \subset \mathbf{G}(d-2, {\proj}V_{d+1})\times|{\Oline}(d)|
\end{equation}
of pairs $({\proj}W, \mathbf{p})$ such that 
$\mathbf{p}=p_1+\cdots+p_d$ is a distinct $d$ point set on ${\proj}^1$ and 
${\proj}W$ is a hyperplane of the $(d-1)$-plane 
${\proj}U_{\mathbf{p}}=\langle \phi_0(p_1), \cdots, \phi_0(p_d)\rangle$. 
We arrived at the birational equivalence 
\begin{equation}
\mathcal{X}_d/{\PGL}_3 \sim \mathcal{Z}_d/{\PGL}_2. 
\end{equation}
Now we prove the 

\begin{proposition}\label{prop 3.3}
If $d\geq5$ is odd, the ${\PGL}_2$-quotient of $\mathcal{Z}_d$ is rational. 
Hence $|L_{1,d}|/{\autF}$ is rational too. 
\end{proposition}

\begin{proof}
The morphism 
\begin{equation}\label{eqn 3.5}
\mathcal{Z}_d \to |{\Oline}(d)|, \quad  ({\proj}W, \mathbf{p})\mapsto \mathbf{p} 
\end{equation}
is dominant with the fiber over a general $\mathbf{p}$ being ${\proj}U_{\mathbf{p}}^{\vee}$. 
The vector space $U_{\mathbf{p}}$ is a subspace of 
$V_{d+1}=H^0({\Oline}(d+1))^{\vee}$. 
Since $d+1$ is even, the bundle ${\Oline}(d+1)$ is ${\PGL}_2$-linearized 
so that the ${\PGL}_2$-action on ${\proj}V_{d+1}$ is induced by a ${\PGL}_2$-representation on $V_{d+1}$. 
Therefore $\mathcal{Z}_d$ is ${\PGL}_2$-isomorphic to 
the projectivization of a ${\PGL}_2$-linearized vector bundle on an open set of $|{\Oline}(d)|$. 
As ${\PGL}_2$ acts almost freely on $|{\Oline}(d)|$, 
the no-name method applied to the fibration \eqref{eqn 3.5} shows that 
\begin{equation}
\mathcal{Z}_d/{\PGL}_2 \sim {\proj}^{d-1}\times(|{\Oline}(d)|/{\PGL}_2). 
\end{equation}
The quotient $|{\Oline}(d)|/{\PGL}_2$ is rational by Katsylo \cite{Ka1}. 
\end{proof}

\noindent
\textit{Proof of Theorem \ref{main'} for} $b\geq4$. 
We may take an odd $d\geq5$ in the range $b \leq 2d+2 \leq 3b+9$. 
By Proposition \ref{prop 2.4}, Lemma \ref{lemma 3.1}, and Lemma \ref{lemma 3.2} we have 
\begin{equation}
|L_{3,b}|/{\autF} \sim {\proj}^{4b+7-2d}\times(|L_{1,d}|/{\autF}). 
\end{equation}
Then $|L_{1,d}|/{\autF}$ is rational by Proposition \ref{prop 3.3}. 
\qed


\section{The case $g\leq7$}\label{sec: g<9}

\subsection{The rationality of $\mathcal{T}_7$}\label{ssec: g=7}

We consider the ${\autF}$-equivariant map $\varphi:|L_{3,3}|\dashrightarrow S^3{\F}_1$ defined in \eqref{map 3}. 
The group ${\autF}$ acts almost transitively on $S^3{\F}_1$, 
with the stabilizer $G$ of a general point $q_1+q_2+q_3$ being isomorphic to $\frak{S}_3$ by the permutation action 
on the set $\{ q_1, q_2, q_3\}$. 
As proved in Lemma \ref{lemma 2.1}, the fiber $\varphi^{-1}(q_1+q_2+q_3)$ is an open set of 
a sub linear system ${\proj}V\subset|L_{3,3}|$. 
Then by the slice method (see \cite{Do}) we have the birational equivalence 
\begin{equation}
|L_{3,3}|/{\autF} \sim  {\proj}V/G. 
\end{equation}
The $G$-action on ${\proj}V$ is induced by a $G$-representation on $V$ 
because the bundle $L_{3,3}$ admits an ${\autF}$-linearization. 
It is well-known that for any linear representation $V'$ of $\frak{S}_3$ the quotient ${\proj}V'/\frak{S}_3$ is rational 
(apply the no-name method for the irreducible decomposition). 
Hence the quotient ${\proj}V/G$ is rational, and Theorem \ref{main'} is proved for $b=3$.

\subsection{The rationality of $\mathcal{T}_5$}\label{ssec: g=5}

We consider the ${\autF}$-equivariant map 
$\varphi_1 : |L_{3,2}| \dashrightarrow S^2({\proj}\mathcal{E})$ defined in \eqref{map 1}. 

\begin{lemma}\label{lemma 4.1}
The group ${\autF}$ acts almost transitively on $S^2({\proj}\mathcal{E})$ 
with the stabilizer $G$ of a general point $\mathbf{q}=(q_{1+}+q_{1-})+(q_{2+}+q_{2-})$ being isomorphic to 
$\frak{S}_2\ltimes(\frak{S}_2\times\frak{S}_2)$. 
\end{lemma}

\begin{proof}
Since ${\autF}$ and  $S^2({\proj}\mathcal{E})$ have the same dimention, it suffices to calculate the stabilizer $G$. 
If $p_{i\pm}\in{\proj}^2$ is the image of $q_{i\pm}$ by the blow-down ${\F}_1\to{\proj}^2$,   
the group $G$ is identified with the group of those $g\in{\PGL}_3$ 
such that for each $i=1, 2$ we have $g(\{ p_{i+}, p_{i-}\}) = \{ p_{j+}, p_{j-}\}$ for some $1\leq j\leq2$. 
\end{proof}

Let $F_i$ be the $\pi$-fiber passing $q_{i\pm}$ and let $p_i=F_i\cap\Sigma$. 
The fiber $\varphi_1^{-1}(\mathbf{q})$ is an open set of the sub linear system 
${\proj}V\subset|L_{3,2}|$ of curves passing $q_{1+}, \cdots, q_{2-}$ and $p_1, p_2$. 
Similarly as Section \ref{ssec: g=7}, the slice method applied to the map $\varphi_1$ implies that 
\begin{equation}
|L_{3,2}|/{\autF} \sim {\proj}V/G, 
\end{equation}
where the $G$-action on ${\proj}V$ is induced by a $G$-representation on $V$. 
Let ${\proj}W\subset{\proj}V$ be the sub linear system defined by 
\begin{equation}
{\proj}W = 2F_1 + 2F_2 + 2\Sigma + |L_{1,0}|. 
\end{equation}
Since the group $G$ preserves the curves $F_1+F_2$ and $\Sigma$, 
the subspace ${\proj}W$ is invariant under the $G$-action. 
Since $G$ is finite, we have a $G$-decomposition 
$V=W\oplus W^{\perp}$ where $W^{\perp}$ is a $G$-invariant subspace. 
The group $G$ acts almost freely on the linear system $|L_{1,0}|$. 
Hence we may apply the no-name lemma for the projection 
${\proj}V\dashrightarrow {\proj}W$ from ${\proj}W^{\perp}$ 
to see that 
\begin{equation}
{\proj}V/G \sim {\C}^9\times({\proj}W/G). 
\end{equation}
The quotient ${\proj}W/G$, being of dimension $2$, is rational by Castelnuovo's theorem. 
This completes the proof of rationality of $\mathcal{T}_5$.


\end{document}